\documentclass[11pt,fleqn,a4paper]{article}
\usepackage{amssymb,latexsym,amsmath,amsfonts,graphicx}
\advance\textheight by \topskip

\def\arg{\textrm{arg}}
\def\th{\textrm{th}}

\def \C{\mathbb{C}}

\def \N{\mathbb{N}}

\def \Z{\mathbb{Z}}

\def \f#1#2#3#4#5#6{
      \: {_{#1}} #2_{#3}\left( \left. #4 \atop #5 \right| #6 \right) }
\def \g#1#2#3#4#5#6{
      \: #1_{#2}^{#3}\left( \left. #4 \atop #5 \right| #6 \right) }

\def \reff#1{(\ref{#1})}
\def \v#1{\vec{#1}}

\numberwithin{equation}{section}
\newtheorem{theorem}{Theorem}[section]
\newtheorem{lemma}[theorem]{Lemma}
\newtheorem{corollary}[theorem]{Corollary}

\newtheorem{Definition}[theorem]{Definition}

\newtheorem{Remark}[theorem]{Remark}
\newenvironment{remark}{\begin{Remark}\rm}{\end{Remark}}
\newtheorem{Example}[theorem]{Example}

\newtheorem{Assumption}[theorem]{Assumption}

\newenvironment{proof}%
{\rm \trivlist \item[\hskip \labelsep{\bf Proof. }]}%
{\hspace*{\fill}$\Box$\endtrivlist}

\begin{document}

    \begin{center} \Large\bf
        Multiple Wilson and Jacobi-Pi\~{n}eiro polynomials
        \footnote{This work was supported by INTAS project 00-272, FWO research project G.0184.02,
        Tournesol project T2002.02, and in part by the Ministry of Science
   and Technology (MCYT) of Spain and the European Regional
   Development Fund (ERDF) through the grant BFM2001-3878-C02-02.}
    \end{center}
    \

    \begin{center}
        \large
        B.\ Beckermann\\
        \normalsize \em
        Laboratoire de Math\'ematiques Appliqu\'ees F.R.E. 2222 (ANO)\\
        UFR Math\'ematiques - M3, Universit\'e de Lille 1\\
        59655 Villeneuve d'Ascq CEDEX\\
        \rm bbecker@ano.univ-lille1.fr\\[3ex]
        \large
        J.\ Coussement\footnote{Research Assistant of the Fund for Scientific Research -- Flanders (Belgium)}\\
        \normalsize \em
        Department of Mathematics, Katholieke Universiteit
        Leuven,\\
        Celestijnenlaan 200 B, 3001 Leuven, Belgium \\
        \rm jonathan.coussement@wis.kuleuven.ac.be\\[3ex]
        \large
        W.\ Van Assche\\
        \normalsize \em
        Department of Mathematics, Katholieke Universiteit
        Leuven,\\
        Celestijnenlaan 200 B, 3001 Leuven, Belgium \\
        \rm walter@wis.kuleuven.ac.be\\[3ex]
    \end{center}

\ \\[1ex]

    \begin{abstract}
    \noindent
    We introduce multiple Wilson polynomials, which give a new example
    of multiple orthogonal polynomials (Hermite-Pad\'e polynomials) of
    type II.  These
    polynomials can be written as a Jacobi-Pi\~{n}eiro transform,
    which is a generalization of the Jacobi transform for Wilson polynomials,
    found by T.H.\ Koornwinder.  Here we need to introduce Jacobi and
    Jacobi-Pi\~{n}eiro polynomials with complex parameters.  Some
    explicit formulas are provided for both
    Jacobi-Pi\~{n}eiro and multiple Wilson polynomials, one of them in terms of
    Kamp\'e de F\'eriet series. Finally we look at some limiting
    relations and construct a part of a multiple AT-Askey table.
    \end{abstract}

\section{Introduction}

Multiple orthogonal polynomials are a generalization of orthogonal
polynomials in the sense that they satisfy orthogonality
conditions with respect to $m\in \N$ measures $\mu_1,\ldots
,\mu_m$ \cite{Nikishin,Apt1,Coussement}.  In this paper $m$ will
always represent the number of weights.  Multiple orthogonal
polynomials arise naturally in the theory of simultaneous rational
approximation, in particular in Hermite-Pad\'e approximation of a
system of $m$ (Markov) functions \cite{Mahler,Bruin1,Bruin2}.

There are two types of multiple orthogonal polynomials.  In the
present paper we only consider multiple orthogonal polynomials of
type II.  Let $\v n = (n_1,n_2,\ldots ,n_m)$ be a vector of $m$
nonnegative integers, which is called a {\em multi-index} with
length $|\v n| := n_1+n_2+\cdots+n_m$.  Furthermore let
$\Gamma_1,\ldots,\Gamma_m$ be the supports of the $m$ measures. A
multiple orthogonal polynomial $p_{\v n}$ of type II with respect
to the multi-index $\v n$, is a (nontrivial) polynomial of degree
$\leq |\v n|$ which satisfies the orthogonality conditions
   \begin{equation}
   \label{system}
   \int_{\Gamma_j} p_{\v n}(t)\: t^k \: d\mu_j (t) =0,
   \qquad k=0,\ldots ,n_j-1,\quad j=1,\ldots ,m.
   \end{equation}
Notice that the measures in (\ref{system}) are not necessarily
supposed to be positive. In case of a complex orthogonality
relation, one usually refers to $p_{\v n}$ as a {\em formal}
multiple orthogonal polynomial.

Equation (\ref{system}) leads to a system of $|\v n|$ homogeneous
linear relations for the $|\v n|+1$ unknown coefficients of $p_{\v
n}$. A basic requirement in the study of such multiple orthogonal
polynomials is that there is (up to a scalar multiplicative
constant) a unique solution of system (\ref{system}). We call $\v
n$ a {\em normal index} for $\mu_1,\ldots,\mu_m$ if any solution
of (\ref{system}) has exact degree $|\v n|$ (which implies
uniqueness). Let $m_k^{(j)}=\int_{\Gamma_j}t^k\: d\mu_j$ be the
$k$th moment of the measure $\mu_j$. Further set
\begin{eqnarray}
\label{momentmatrix} D_{\v n} & = & \left(
\begin{array}{c}
D_1(n_1)\\ \vdots\\ D_r(n_r)
\end{array} \right),
\end{eqnarray}
where
\begin{displaymath}
D_j(n_j) = \left(
\begin{array}{cccc}
m_0^{(j)} & m_1^{(j)} & \ldots & m_{|\vec{n}|-1}^{(j)}\\ \vdots &
\vdots & & \vdots\\ m_{n_j-1}^{(j)} & m_{n_j}^{(j)} & \ldots &
m_{|\vec{n}|+n_j-2}^{(j)}
\end{array} \right)
\end{displaymath}
is an $n_j\times (|\vec{n}|)$ matrix of moments of the measure
$\mu_j$.  Then $D_{\v n}$ is the matrix of the linear system
\reff{system} without the last column.  It is known and easily
verified that the multi-index $\vec{n}=(n_1,\ldots,n_r)$ is normal
if and only if this matrix has rank $|\vec{n}|$.  When every
multi-index is normal we call the system of measures a {\em
perfect system}. For perfect systems, the multiple orthogonal
polynomials of type II satisfy a recurrence relation of order
$m+1$. The proof is similar as in the case of orthogonal
polynomials, see for instance \cite{Apt1}. Because of this
recurrence relation, formal multiple orthogonal polynomials are a
useful tool in the spectral theory of non-symmetric linear
difference operators.

In the literature one can find some examples of multiple
orthogonal polynomials with respect to positive measures on the
real line which have the same flavor as the classical orthogonal
polynomials. Two classes of measures have been analyzed in more
detail and are known to form a perfect system, see for instance
the monograph \cite{Nikishin} or the survey given in
\cite{Coussement}. The first class consists of Angelesco systems
where the supports of the measures are disjoint intervals. In the
second class of so-called AT systems, the supports of the $m$
measures coincide, and the functions $d\mu_j(x)/d\mu_1(x)$ for
$j=1,\ldots,m$ form an algebraic Chebyshev system
\cite[Section~IV.4]{Nikishin} on the convex hull of the support.
In the continuous case (where the measures can be written as
$d\mu_j(x)=w_j(x)dx$, with $w_j$ the weight function of the
measure $\mu_j$) there are multiple Hermite, multiple Laguerre I
and II, Jacobi-Pi\~{n}eiro, multiple Bessel, Jacobi-Angeleso,
Jacobi-Laguerre and Laguerre-Hermite polynomials, see
\cite{Nikishin,Coussement,Apt} and the references therein. Some
classical discrete examples are multiple Charlier, multiple
Kravchuk, multiple Meixner I and II and multiple Hahn
\cite{Arvesu}. For all these examples there exists a first-order
{\em raising operator}, based on the differential operator $D$ or
the difference operators $\Delta$ and $\nabla$, and a {\em
Rodrigues formula} so that they can be called classical. Moreover,
there exist differential or difference equations of order $m+1$
(with polynomial coefficients)\cite{Apt}. The recurrence relations
of order $m+1$ are known explicitly  for these examples in the
case $m\le 2$. Finally we mention that there also exist some
examples of multiple orthogonal polynomials associated with
modified Bessel functions \cite{Couss1,Couss2,Yak} which can be
called classical.

In Subsection \ref{defjacpin} we recall the definition of one of
these examples, namely Jacobi-Pi\~{n}eiro polynomials $P^{(\v
\alpha,\beta)}_{\v n}$, which are orthogonal with respect to the
weights $w_j(x)=x^{\alpha_j}(1-x)^\beta$ on $[0,1]$,
$\alpha_j,\beta\geq -1$. These polynomials reduce to the classical
Jacobi polynomials (shifted to the interval $[0,1]$) for $m=1$. We
show in Subsection~\ref{Jacobi} that Jacobi polynomials remain
formal orthogonal polynomials for complex parameters
$\alpha_1,\beta$, the corresponding complex orthogonality relation
being obtained via an analytic extension of the Beta function in
both variables. As we show in Subsection~\ref{defjacpin}, also
Jacobi-Pi\~{n}eiro polynomials with complex parameters are formal
multiple orthogonal polynomials of type II.

In Subsection \ref{multwilpolsection} we then introduce the formal
multiple Wilson polynomials $p_{\v n}(\cdot;a,\v b,c,d)$, which
give a new example of formal multiple orthogonal polynomials of
type II. They are an extension of the formal Wilson polynomials
$p_n(\cdot;a,b,c,d)$ \cite{wils1,wils2} for which we recall the
definition in Subsection \ref{defWilsonsection}. We also mention
that, with some conditions on the complex parameters $a,b,c,d$, we
find the Wilson and Racah polynomials on the top of the Askey
scheme which have real orthogonality conditions.

The formal multiple Wilson polynomials satisfy a complex
orthogonality conditions with res\-pect to $m$ Wilson weights
    \begin{equation}
    \label{Wilson_weight}
        w(t;a,b_j,c,d)= \frac{\Gamma(a+t)\Gamma(a-t)\Gamma(b_j+t)\Gamma(b_j-t)
        \Gamma(c+t)\Gamma(c-t)\Gamma(d+t)\Gamma(d-t)}{\Gamma(2t)\Gamma(-2t)},
    \end{equation}
$j=1,\ldots,m$, where we integrate over the imaginary axis
deformed so as to separate the increasing sequences of poles of
these weight functions from the decreasing ones. Note that there
are some additional conditions on the complex parameters in order
to have that these Wilson weights have only simple poles. We prove
in Theorem~\ref{generaltransform} that the weight functions
\reff{Wilson_weight} form a perfect system if $b_l-b_k\notin\Z$
whenever $l\not= k$. In the same Theorem we show that, for
$\Re(c+d)>0$ and $0<|\Re(t)|<\Re(a)$, the formal multiple Wilson
polynomials can be written as the Jacobi-Pi\~{n}eiro transform
        \begin{equation}
        \label{jacpintrans}
         p_{\v n}(t^2; a,\v b,c,d) =
         k_{\v n}\int_{0}^{1}
          P^{(\v \alpha,\beta)}_{\v n}(u)
          w^{(a-1,\beta)}(u)  K(u,t;a,0,c,d)\: du,
        \end{equation}
where $\v \alpha=(a+b_1-1,...,a+b_m-1)$ and $\beta=c+d-1$.  Here
$k_{\v n}$ is a normalizing constant, $w^{(\alpha,\beta)}(u)=
u^\alpha(1-u)^\beta$ the Jacobi weight and
        \begin{eqnarray}
        \label{kern1}
        K(u,t;a,b,c,d) &=& \frac{ u^{-b-t} }
        {\Gamma(a-t) \Gamma(a+t) \Gamma(c+d)}
        \ \f{2}{F}{1}{c-t,d-t}{c+d}{1-u}
        \end{eqnarray}
some kernel function, independent of $n$. In the scalar case
($m=1$) this formula reduces to a Jacobi transform for the Wilson
polynomials which was already found by T.H.\ Koornwinder
\cite[(3.3)]{Koorn}.  We recall this Jacobi transform in
Subsection \ref{wiljacscalartrans} and give a short proof. Note
that the parameters $a,b_1,\ldots,b_m,c,d$ can take complex
values.

The Jacobi-Pi\~{n}eiro transform \reff{jacpintrans} is the key
formula of this paper. In Subsection \ref{defjacpin} we obtain two
new hypergeometric representations for the Jacobi-Pi\~{n}eiro
polynomials starting from the Rodrigues formula.   Applying the
Jacobi-Pi\~{n}eiro transform \reff{jacpintrans} we then also find
two explicit formulas for the formal multiple Wilson polynomials
(see Subsection \ref{multwilpolsection}). One of them is in terms
of Kamp\'e de F\'eriet series \cite{Srivastava}.

In Section \ref{multAskeysection} we only consider the cases where
we obtain real orthogonality conditions, namely multiple Wilson
and multiple Racah. By some limit relations we then recover
hypergeometric re\-presentations for the examples of multiple
orthogonal polynomials of type II, mentioned above.  We also
introduce some new examples like multiple dual Hahn,  multiple
continuous dual Hahn, and multiple Meixner-Pollaczek.  As a
result, we finally construct a (still incomplete) multiple
AT-Askey table similar to the well known Askey scheme for
classical orthogonal polynomials.



\section{Jacobi and Wilson polynomials}

\subsection{Formal (shifted) Jacobi polynomials}
\label{Jacobi}

The (shifted) Jacobi polynomials $P_n^{(\alpha,\beta)}$ are a
classical example of continuous orthogonal polynomials.  Suppose
$\alpha,\beta > -1$, then these polynomials are orthogonal with
respect to the Jacobi weight function $w^{(\alpha,\beta)}(x) =
x^\alpha (1-x)^\beta$ on the interval $[0,1]$.  These polynomials
have the explicit expressions \cite{Chihara}
    \begin{eqnarray}
        \label{onehyper} P_n^{(\alpha,\beta)}(x)&=& \frac{(\alpha+1)_n}{n!}
        \f{2}{F}{1}{-n,\alpha+\beta+n+1}{\alpha+1}{x},\\
        \label{twohyper} &=& \frac{(\alpha+1)_n}{n!}\: (1-x)^{-\beta}
        \f{2}{F}{1}{\alpha+1+n,-\beta-n}{\alpha+1}{x},
    \end{eqnarray}
where the second expression is obtained by Euler's formula
\cite[15.3.3]{Abra},\cite{Erdelyi}. We claim that the (shifted)
Jacobi polynomials $P_n^{(\alpha,\beta)}$ are still formal
orthogonal polynomials if we allow complex parameters
$\alpha,\beta,\alpha+\beta+1 \in \C \setminus \{-1,-2,\ldots\}$ in
formula \reff{onehyper}.  This was already mentioned in
\cite[Theorem 2.1]{Arno}, but we give a different proof.

\begin{figure}[t]
\begin{center}
\includegraphics[scale=0.5]{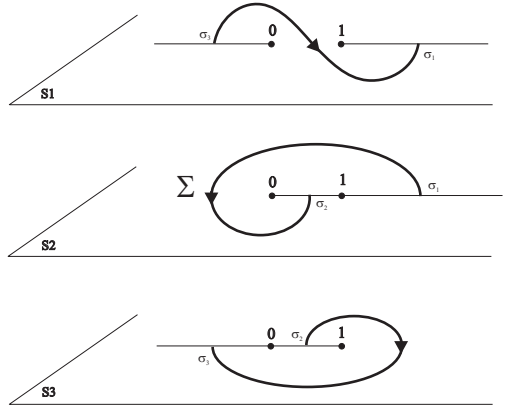}
\end{center}
\caption{\label{contour}{\em The contour $\Sigma$ on the
appropriate Riemann surface for the function
$t^{z-1}(1-t)^{w-1}$.}}
\end{figure}

In order to prove this claim, we require an integral
representation for the meromorphic continuation in both variables
of the beta function. Recall, e.g., from \cite[\S~1.1]{Erdelyi}
that $\Gamma$ is meromorphic in $\mathbb C$, with simple poles at
$0,-1,-2,...$, and hence the beta function \cite[\S~1.5]{Erdelyi}
    \[B(z,w)=\frac{\Gamma(z)\Gamma(w)}{\Gamma(z+w)}\]
is meromorphic both in $z$ and $w$, with simple poles at
$z,w=0,-1,-2,...$. From \cite[\S~1.5]{Erdelyi} we have the
integral representation
    \begin{equation}
    \label{betadef}
    B(z,w)=\int_0^1 t^{z-1}(1-t)^{w-1}dt,\qquad \Re (z)>0,\Re (w) >
    0.
    \end{equation}
In order to obtain a representation valid for general $z,w \in
\C\setminus \Z$ (compare with the Pochhammer formula
\cite[1.6.(7)]{Erdelyi}), we follow \cite[\S\ 3.4]{Hochstadt} and
consider three sheets $S1$, $S2$ and $S3$ of the appropriate
Riemann surface for the function $t^{z-1}(1-t)^{w-1}$ (in the
variable $t$) so that
    \[
    \left\{ \begin{array}{lll}
    -\pi < \arg (t)<\pi,\quad -\pi < \arg (1-t)<\pi, & \qquad & \mbox{for\ }t\in
    S1\setminus \left\{(-\infty,0] \cup [1,+\infty)\right\},\\
    \ \ 0 < \arg (t)<2\pi,\quad \ \pi < \arg (1-t)<3\pi, & \qquad & \mbox{for\ }t\in
    S2\setminus \left\{[0,+\infty) \cup [1,+\infty)\right\},\\
    \ \ \pi < \arg (t)<3\pi,\quad \ 0 < \arg (1-t)<2\pi, & \qquad & \mbox{for\ }t\in
    S3\setminus \left\{(-\infty,0] \cup (-\infty,1]\right\}.
    \end{array} \right.
    \]
Furthermore we choose a closed contour $\Sigma$ as in Figure
\ref{contour} where $\sigma_1,\sigma_2,\sigma_3$ are the
transition points between the three sheets.  Note that the
function $t^{z-1}(1-t)^{w-1}$ is analytic on $\Sigma$. For the
beta function we then have
    \begin{equation}
    \label{beta}
    B(z,w)=(1-e^{2\pi iz})^{-1}(1-e^{2\pi
    iw})^{-1}\int_{\Sigma}t^{z-1}(1-t)^{w-1}dt, \quad z,w\in \C \setminus
    \Z.
    \end{equation}
Indeed, if $\Re (z)>0,\Re (w)>0$ then the path of integration in
(\ref{beta}) can be shrinked in order to approach the interval
$[0,1]$, leading to formula (\ref{betadef}). In particular, if $z$
and/or $w$ is a strictly positive integer, we can obtain $B(z,w)$
by taking limits in (\ref{beta}).

Now we prove our claim that (shifted) Jacobi polynomials with
complex parameters are formal orthogonal polynomials.

    \begin{theorem}
    \label{formal_Jacobi}
    Let $\alpha,\beta,\alpha+\beta+1 \in \C \setminus
    \{-1,-2,\ldots\}$, and consider the (complex) measure
    $\mu^{(\alpha,\beta)}$ defined by
        \[
           \int h(x) \, d\mu^{(\alpha,\beta)}(x) =
           \lim_{z\to\alpha,w\to\beta}(1-e^{2\pi iz})^{-1}(1-e^{2\pi iw})^{-1}
           \int_{\Sigma} h(t) t^z (1-t)^w dt .
        \]
    Then $\mu^{(\alpha,\beta)}$ forms a perfect system, with the
    corresponding $n$th formal orthogonal polynomial given by the
    (shifted) Jacobi polynomial $P_n^{(\alpha,\beta)}$.
    \end{theorem}
    \begin{proof}
    From (\ref{beta}) we obtain for the $k$th moment
       \[
       \int x^k \, d\mu^{(\alpha,\beta)}(x) =
       B(\alpha+1+k,\beta+1).
       \]
    The restrictions on the parameters $\alpha$ and $\beta$ guarantee that
    the determinant of the moment matrix,
        \[
        \det \Bigl(B(\alpha+r+s-1,\beta+1)\Bigr)_{1\le r,s\le
        n}=\prod_{i=1}^n\frac{\Gamma(\alpha+i)\Gamma(\beta+i)}
        {\Gamma(\alpha+\beta+n+i)}\prod_{1\le r<s\le n}(s-r),
        \]
    is different from zero. A proof for this formula uses
    Theorem 1 and 2 in \cite{Corput1,Corput2}  (the moment matrix in
    the non-shifted case can be found in \cite[6.71.5]{Szego}).
    Thus $\mu^{(\alpha,\beta)}$ forms a perfect system. Moreover,
    for $0\le j\le n-1$ we have according to \reff{onehyper} and
    \reff{beta}
        \begin{eqnarray*}
        \lefteqn{
        \int (1-x)^j
        P_n^{(\alpha,\beta)}(x) \, d\mu^{(\alpha,\beta)}(x)
        }\\
        & = &
        \frac{(\alpha+1)_n}{n!}\sum_{k=0}^n \frac{(-n)_k(\alpha+\beta+n+1)_k}{(\alpha+1)_k k!}
        \ B(\alpha+k+1,\beta+j+1)\\
        & = &
        \frac{(-1)^n\Gamma(\alpha+n+1)\Gamma(\beta+j+1)}{\Gamma(\alpha+\beta+n+1)}\sum_{k=0}^n
        \frac{(-1)^{n-k}}{ k!(n-k)!
        }\prod_{\ell=j+2}^{n}(\alpha+\beta+k+\ell).
        \end{eqnarray*}
    The sum on the right hand side is the divided difference
    $g_{j,n}[0,1,...,n]$ of the polynomial
    $g_{j,n}(x)=\prod_{l=j+2}^{n}(\alpha+\beta+x+l)$ of degree $n-j-1$.
    Since $n-j-1<n$, this is equal to $0$, which proves that the
    Jacobi polynomial is a formal orthogonal
    polynomial with respect to $\mu^{(\alpha,\beta)}$.
    \end{proof}

\subsection{Formal Wilson polynomials}
\label{defWilsonsection}

In \cite{wils1} J.A.\ Wilson introduced the (formal) Wilson
polynomials
    \begin{equation}
    \label{Wilson_definition}
       p_n(t^2; a,b,c,d) = (a+b)_n (a+c)_n (a+d)_n\
       \f{4}{F}{3}{-n,a+b+c+d+n-1,a-t,a+t}{a+b,a+c,a+d}{1}.
    \end{equation}
Due to Whipple's identities \cite[Theorem 3.3.3]{Andrews}, one can
show that these formal Wilson polynomials are symmetric in all
four complex parameters $a,b,c,d$. Furthermore, with some
conditions on these four parameters, the polynomials satisfy a
complex orthogonality with respect to the Wilson
 weight function
    \begin{equation*}
    \nonumber
        w(t;a,b,c,d)= \frac{\Gamma(a+t)\Gamma(a-t)\Gamma(b+t)\Gamma(b-t)
        \Gamma(c+t)\Gamma(c-t)\Gamma(d+t)\Gamma(d-t)}{\Gamma(2t)\Gamma(-2t)}.
    \end{equation*}
Suppose that
    \begin{equation}
    \label{conditions}
    2a,a+b,a+c,a+d,2b,b+c,b+d,2c,c+d,2d \notin \{0,-1,-2,\ldots \},
    \end{equation}
so that the Wilson weight has only simple poles, and that
    \begin{equation}
    \label{condition2}
    a+b+c+d\notin \{0,-1,-2,\ldots \}.
    \end{equation}
Furthermore denote by the contour ${\cal C}$ the imaginary axis
deformed so as to separate the increasing sequences of poles
($\{a+k\}_{k=0}^\infty,\{b+k\}_{k=0}^\infty,\{c+k\}_{k=0}^\infty,\{d+k\}_{k=0}^\infty$)
from the decreasing ones
($\{-a-k\}_{k=0}^\infty,\{-b-k\}_{k=0}^\infty,\{-c-k\}_{k=0}^\infty,
\{-d-k\}_{k=0}^\infty$), and define the Wilson measure
$\mu^{(a,b,c,d)}$ by
    \begin{equation}
    \int h(t) \, d\mu^{(a,b,c,d)}(t) = \int_{\mathcal C} h(t^2)
    w(t;a,b,c,d) \, dt .
    \end{equation}
Wilson \cite{wils1} shows the complex orthogonality relations
    \[
       \int p_m(t;a,b,c,d)p_n(t;a,b,c,d) \, d\mu^{(a,b,c,d)}(t)
       = \delta_{mn} 2 i M_n,
    \]
where
    \begin{eqnarray*}
    M_n & = & 2\pi n!\ (a+b+c+d+n-1)_n\\
    & & \quad \times \frac{\Gamma(a+b+n)\Gamma(a+c+n)
    \Gamma(a+d+n)\Gamma(b+c+n)\Gamma(b+d+n)\Gamma(c+d+n)}{\Gamma(a+b+c+d+2n)}.
    \end{eqnarray*}
The perfectness of the singleton system $\mu^{(a,b,c,d)}$ follows
from Lemma~\ref{lemmaperfect} below.

In some cases we obtain real orthogonality conditions with respect
to positive measures on the real line \cite{wils1}. If
$\Re(a),\Re(b),\Re(c),\Re(d)>0$ and $a,b,c$ and $d$ are real
except possibly for conjugate pairs, then $\cal C$ can be taken to
be the imaginary axis and we obtain the real orthogonality
    \[
    \int_0^{\infty} p_m(-x^2;a,b,c,d)p_n(-x^2;a,b,c,d)
    \left|\frac{\Gamma(a+ix)\Gamma(b+ix)\Gamma(c+ix)\Gamma(d+ix)}{\Gamma(2ix)}\right|^2
    dx = \delta_{mn} M_n.
    \]
Another case is when $a<0$ and $a+b,a+c,a+d$ are positive or a
pair of complex conjugates occurs with positive real parts. We
then get the same positive continuous weight function where some
positive point masses are added.
 In these cases we obtain the Wilson polynomials
$W_n(x^2;a,b,c,d):=p_n(-x^2;a,b,c,d)$ (see, e.g., \cite{Koekoek}),
which are real for real $x^2$.

Finally in the case that one of $a+b,a+c,a+d$ is equal to
$-N+\epsilon$, with $N$ a nonnegative integer, we obtain a purely
discrete orthogonality after dividing by $\Gamma(-N+\epsilon)$ and
letting $\epsilon\to 0$. If we take the substitution $t\to x+a$
and the change of variables $\alpha=a+b-1, \beta=c+d-1,
\gamma=a+d-1, \delta=a-d$ we then obtain the Racah polynomials
    \[
        R_n(\lambda(x); \alpha,\beta,\gamma,\delta) =
       \f{4}{F}{3}{-n,n+\alpha+\beta+1,-x,x+\gamma+\delta+1}{\alpha+1,\beta+\delta+1,
       \gamma+1}{1},
    \]
where $\lambda(x)=x(x+\gamma+\delta+1)$ and $\alpha +1=-N$ or
$\beta+\delta+1=-N$ or $\gamma+1=-N$.  They satisfy the discrete
orthogonality
    \[
    \sum_{x=0}^N\frac{(\alpha+1)_x(\gamma+1)_x(\beta+\delta+1)_x(\gamma+\delta+1)_x((\gamma+\delta+3)/2)_x}
    {(-\alpha+\gamma+\delta+1)_x(-\beta+\gamma+1)_x((\gamma+\delta+1)/2)_x(\delta+1)_xx!}\
    R_n(\lambda(x); \alpha,\beta,\gamma,\delta)\:(\lambda(x))^k=0,
    \]
$k=0,\ldots,n-1$.  Necessary and sufficient conditions for the
positivity of the weights are quite messy.  An example of
sufficient conditions is given in \cite[(3.5)]{wils1}.

The formal Wilson polynomials contain as limiting cases several
families of orthogonal polynomials \cite[Section
4]{wils1}\cite[Chapter 2]{Koekoek} like Hahn, dual Hahn, Meixner,
Krawtchouk, Charlier, continuous Hahn, continuous dual Hahn,
Meixner-Pollaczek, Jacobi, Laguerre and Hermite polynomials.
Together they form the Askey scheme of hypergeometric orthogonal
polynomials. Moreover, there exist $q$-analogues
\cite{Koekoek} for all these polynomials.  The results of this
paper should probably extend to this more general case.


\subsection{Formal Wilson polynomials as a Jacobi transform}
\label{wiljacscalartrans}

We now recall an integral relation between the Jacobi and the
formal Wilson polynomials.  We also give a short proof which will
help us to find explicit expressions for the formal multiple
Wilson polynomials in the next section.

    \begin{theorem}[T.H.\ Koornwinder]
    Suppose that the conditions \reff{conditions} and \reff{condition2}
    hold, and that
    $\Re(c+d)>0$, $0<|\Re(t)|<\Re(a)$. Then we have
        \begin{equation}
        \label{scal_trans}
        p_n(t^2; a,b,c,d) = k_n
        \int_{0}^{1}
        P^{(\alpha,\beta)}_n(u)  w^{(\alpha,\beta)}(u) K(u,t;a,b,c,d)\ du,
        \end{equation}
    with $\alpha=a+b-1$ and $\beta=c+d-1$, the Jacobi weight
    function $w^{(\alpha,\beta)}(u)=u^\alpha(1-u)^\beta$, the
    constant $k_n = n!\ \Gamma(a+c+n) \Gamma(a+d+n)$
    and the kernel \reff{kern1}.
    \end{theorem}

    \begin{remark}
    T.H.\ Koornwinder already mentioned this Jacobi transform for the
    Wilson polynomials in \cite[(3.3)]{Koorn}.  In order to
    see that these formulas coincide, one has to
    make some parameter changes and substitute $1-\th^2s\to u$.
    Notice also that \reff{scal_trans} is a particular case of a
    formula due to Meijer \cite[p.103]{meij}.
    \end{remark}

    \begin{proof}
    By comparing the explicit formulas \reff{onehyper} and
    \reff{Wilson_definition} we see that it is sufficient to prove
    that, for $\ell\in \N$,
        \begin{equation}
        \label{bewijstransform}
        \int_{0}^{1} u^\ell w^{(a-1,c+d-1)}(u) K(u,t;a,0,c,d)\: du
          = \frac{(a-t)_\ell (a+t)_\ell}{ \Gamma(a+c+\ell)  \Gamma(a+d+\ell)}.
        \end{equation}
    By definition of the kernel \reff{kern1} we have
        \[ K(u,t;a,0,c,d) 
        = (a-t)_\ell
        (a+t)_\ell \ K(u,t;a+\ell,0,c,d) , \]
    and $u^\ell w^{(a-1,c+d-1)}(u) = w^{(a+\ell-1,c+d-1)}(u)$.
    Furthermore, Euler's formula gives the symmetry
    $K(u,-t;a+\ell,0,c,d)=K(u,t;a+\ell,0,c,d)$.
    By replacing $a+\ell$ by $a$ and possibly $t$ by $-t$, we see that it
    remains to prove that, for $0<\Re(t)<\Re(a)$ and $\Re(c+d)>0$,
        \begin{equation} \label{remains}
         \int_{0}^{1} w^{(a-1,c+d-1)}(u) K(u,t;a,0,c,d)\: du =
        \frac{1}{\Gamma(a+c)  \Gamma(a+d)} .
        \end{equation}
    Denoting the left hand side of \reff{remains} by $J$, we have by
    definition of the kernel
        \[J=
        \frac{1}{\Gamma(a-t) \Gamma(a+t) \Gamma(c+d)}
        \int_{0}^1w^{(a-t-1,c+d-1)}(u)
        \f{2}{F}{1}{c-t,d-t}{c+d}{1-u} du.
        \]
    In order to exchange the order of summation and integration,
    we notice that, if $\Re(C-A-B) >0$,
        \[\lim_{y \to 1-} \ \max_{v \in [0,1]}
          \left| \f{2}{F}{1}{A,B}{C}{y v} - \f{2}{F}{1}{A,B}{C}{v}
          \right| = 0 ,
        \]
    which follows by observing that
    $\frac{\Gamma(A+k)\Gamma(B+k)}{\Gamma(C+k)\Gamma(1+k)}=k^{A+B-C-1}(1+\mathcal
    O(k^{-1}))$, $_{k\to \infty}$, by Stirling's formula. Consequently,
    using the assumptions $\Re(2t)>0$, $\Re(a-t)>0$, $\Re(c+d)>0$
    together with \reff{betadef} we obtain by uniform convergence
        \begin{eqnarray*}
        J &=&
        \frac{1}{\Gamma(a-t) \Gamma(a+t) \Gamma(c+d)}
        \lim_{y \to 1-} \int_{0}^1w^{(a-t-1,c+d-1)}(u)
            \f{2}{F}{1}{c-t,d-t}{c+d}{y(1-u)} du\\
        &=&
        \lim_{y \to 1-} \sum_{k=0}^\infty \frac{ (c-t)_k  (d-t)_k\  y^k}
                { \Gamma(a-t)  \Gamma(a+t)  \Gamma(c+d+k) k!}
                \int_{0}^1  w^{(a-t-1,c+d+k-1)}(u)\: du\\
        &=&
        \lim_{y \to 1-} \frac{1}{ \Gamma(a+t)  \Gamma(a+c+d-t)}
            \f{2}{F}{1}{c-t,d-t}{a+c+d-t}{y} .
        \end{eqnarray*}
    By assumption on the parameters, $\Re(a+c+d-t-(c-t+d-t))>0$,
    and hence we get from the Gauss formula
    \cite[Theorem 2.2.2]{Andrews},\cite[15.1.20]{Abra}
        \[
          J = \frac{1}{ \Gamma(a+t)  \Gamma(a+c+d-t)}
        \f{2}{F}{1}{c-t,d-t}{a+c+d-t}{1} =
        \frac{1}{\Gamma(a+d)\Gamma(a+c)},
        \]
    as claimed in \reff{remains}.
\end{proof}


\section{Multiple Wilson as a Jacobi-Pi\~{n}eiro transform}
\label{multwilsection}

\subsection{Jacobi-Pi\~{n}eiro with complex parameters}
\label{defjacpin}
The Jacobi-Pi\~{n}eiro polynomials are defined by a Rodrigues
formula
    \begin{equation}
    \label{rodriguesjacpin}
    P_{\v n}^{(\v \alpha,\beta)}(t)  =
    \frac{1}{\v n!}(1-t)^{-\beta}
       \prod_{j=1}^m \left( t^{-\alpha_j} \frac{d^{n_j}}
       {dt^{n_j}} t^{n_j+\alpha_j}\right) (1-t)^{|\v
       n|+\beta},
    \end{equation}
where $\v n!=\prod_{j=1}^m n_j!$. It is well known \cite{Pineiro}
that, provided that $\alpha_j>-1,\beta>-1$ and $\alpha_\ell -
\alpha_k \not \in \Z$ for $\ell\neq k$, the Jacobi-Pi\~{n}eiro
polynomials are orthogonal with respect to the (positive) Jacobi
weights $w^{(\alpha_j,\beta)}$, $j=1,\ldots ,m$, on the interval
$[0,1]$. Notice that these weights form an AT system, and hence we
obtain a perfect system of measures. Similar to
Theorem~\ref{formal_Jacobi}, we show below that for complex
parameters we keep formal multiple orthogonal polynomials of type
II. Here we use the measures $\mu^{(\alpha_j,\beta)}$ of
Theorem~\ref{formal_Jacobi} which have as support the contour
$\Sigma$, but can be reduced to complex Jacobi weights
$w^{(\alpha_j,\beta)}$ on the interval $[0,1]$ in the case
$\Re(\alpha_j)>-1,\Re(\beta)>-1$.

    \begin{theorem}
    \label{perfect_Jacobi_Pinero}
    Let $\alpha_j,\beta,\alpha_j+\beta\in \mathbb C \setminus \{
    -1,-2,...\}$. Then the measures
    $\mu^{(\alpha_j,\beta)}$, $j=1,...,m$, of
    Theorem~\ref{formal_Jacobi} form a perfect system. The corresponding
    formal multiple orthogonal polynomial of type II with respect
    to the multi-index $\v n$ is given by \reff{rodriguesjacpin}.
    \end{theorem}

    \begin{proof}
    From Theorem 1 and 2 in \cite{Corput1,Corput2}
    we obtain for the
    determinant of the matrix of moments \reff{momentmatrix} the expression
        \begin{eqnarray*}
        D_{\v n}^{\v \alpha,\beta}
        & = & \left(\prod_{i=1}^{|\v n|}\Gamma(\beta+i)\right)
        \left(\prod_{j=1}^{m}\prod_{i=1}^{n_j}\frac{\Gamma(\alpha_j+i)}{\Gamma(\alpha_j+\beta+\v
        n+i)}\right)\\
        & & \qquad \times \left(\prod_{j=1}^{m}\prod_{1\le r<s\le n_j}(s-r)\right)
        \left(\prod_{1\le i<j\le
        m}\prod_{s=1}^{n_i}\prod_{r=1}^{n_j}(\alpha_j-\alpha_i+r-s\right).
        \end{eqnarray*}
    For our choice of parameters, this expression is different from
    zero, and hence every multi-index is normal.

    Our claim on the (formal) orthogonality of Jacobi-Pi\~{n}eiro polynomials
    will be shown by induction on $m$. For $m=1$, equation
    \reff{rodriguesjacpin} reduces to \reff{onehyper}, and
    the claim follows from Theorem~\ref{formal_Jacobi}. For $m>1$,
    we observe first that, by the Rodrigues formula \reff{rodriguesjacpin},
        \begin{equation}
        \label{hulpje}
        P_{\v n}^{(\v \alpha,\beta)}(t)  =
        \frac{t^{-\alpha_1}(1-t)^{-\beta}}{n_1!}
        \frac{d^{n_1}}{dt^{n_1}}  \left(
        t^{\alpha_1+n_1}(1-t)^{\beta+n_1}
            P_{(n_2,...,n_m)}^{((\alpha_2,...,\alpha_m),\beta+n_1)}(t)
            \right).
        \end{equation}
    By induction hypothesis,
    $P_{(n_2,...,n_m)}^{((\alpha_2,...,\alpha_m),\beta+n_1)}$ is
    a polynomial of degree $|\v n|-n_1$, and thus there exist
    scalars $c_j$ with
        \[
        P_{(n_2,...,n_m)}^{((\alpha_2,...,\alpha_m),\beta+n_1)}(t)
        =\sum_{j=0}^{|\v n|-n_1} c_j
        P^{(\alpha_1+n_1,\beta+n_1)}_j(t) .
        \]
    From the Rodrigues formula for $m=1$ and \reff{hulpje} we conclude that
        \[
        P_{\v n}^{(\v \alpha,\beta)}(t)  =
        \sum_{j=0}^{|\v n|-n_1} c_j
        P^{(\alpha_1,\beta)}_{j+n_1}(t),
        \]
    implying that $P_{\v n}^{(\v \alpha,\beta)}$ is a polynomial of
    degree $|\v n|$ with $\int t^k P_{\v n}^{(\v \alpha,\beta)}(t)
    \, d\mu^{(\alpha_1,\beta)}=0$ for $k=0,1,...,n_1-1$ by
    Theorem~\ref{formal_Jacobi}. The other orthogonality
    conditions are obtained by obser\-ving that \reff{rodriguesjacpin}
    remains invariant if one changes the order in the product in
    \reff{rodriguesjacpin}.
    \end{proof}

With help of the Leibniz rule applied to the Rodrigues formula,
the authors in \cite{Coussement} derive for $m=2$ an explicit
expression in terms of $m=2$ sums. We now give a generalization of
this formula for $m \geq 2$, using the notation
    \begin{eqnarray}
    \label{nieuweformule}
    \lefteqn{
    \g{\mathcal{M}}{q,\v n}{p;r}{\v f ;\v g_1:\cdots :\v g_r}
    {\v \phi ; \v \psi_1:\cdots :\v \psi_r}{\v x}
    } & & \\
    & := & \nonumber
    \underbrace{\sum_{k_1=0}^{n_1}\cdots\sum_{k_m=0}^{n_m}}_{\mbox{$m$ sums}}
    \frac{\prod\limits_{\ell=1}^{p}(f_{\ell})_{|\v k|}}
    {\prod\limits_{\ell=1}^{q}(\phi_{\ell})_{|\v k|}}
    \frac{\prod\limits_{i=1}^{r}(g_{i,1})_{|\v k|-k_1}
    \cdots (g_{i,m-1})_{k_m}}
    {\prod\limits_{i=1}^{r}(\psi_{i,1})_{|\v k|-k_1}
    \cdots (\psi_{i,m-1})_{k_m}}
    \prod_{j=1}^m(-n_j)_{k_j}\frac{x_j^{k_j}}{k_j!},
    \end{eqnarray}
where $\v k=(k_1,\ldots,k_m)$, $\v n\in
\left(\N\cup\{0\}\right)^m$, $\v f\in\C^p$, $\v \phi\in\C^q$ and
$\v g_1,\ldots,\v g_r,\v \psi_1,\ldots,\v \psi_r\in\C^{m-1}$. We
also give in (\ref{hypereerstejac}) another new
 explicit expression for the formal
Jacobi-Pi\~{n}eiro polynomials which reduces to \reff{twohyper} if
$m=1$.

    \begin{theorem}
    \label{expressionsJacPin}
    Let $\v e=(1,\ldots,1)$ be a multi-index of length $m$ and $s(\v n)=(n_1,n_1+n_2,\ldots,|\v n|)$.
    Denote by $\v v^{(j)}$ the vector $\v v$
    without the $j$th component.  For the
    Jacobi-Pi\~{n}eiro polynomials we have the hypergeometric
    representations
        \begin{equation}
        \label{hypertweedejac}
        P_{\v n}^{(\v \alpha,\beta)}(t)
        =
        \frac{(\v \alpha+\v e)_{\v n}}{\v n!}
        \g{\mathcal{M}}{1,\v n}{1;2}{(\alpha_1+\beta+n_1+1);(\v \alpha+\v n+\v e)^{(m)}:
        (\v \alpha+s(\v n)+ (\beta+1)\v e)^{(1)}}
        {(\alpha_1+1) ;(\v \alpha+\v e)^{(1)}:(\v \alpha+s(\v n)+(\beta+1)\v e)^{(m)}}{t\v e}
        \end{equation}
    and
        \begin{equation}
        \label{hypereerstejac}
        P_{\v n}^{(\v \alpha,\beta)}(t)
         =
        \frac{(\v \alpha+\v e)_{\v n}}{\v n!} \ (1-t)^{-\beta}
        \f{m+1}{F}{m}{\v \alpha+\v n + \v e,-\beta-|\v n|}{\v \alpha + \v e}{t} ,
        \end{equation}
    where $\v n!=\prod_{j=1}^m n_j!$ and
    $(\v \alpha+\v e)_{\v n} = \prod_{j=1}^m(\alpha_j+1)_{n_j}$.
    \end{theorem}
    \begin{proof}
    We prove \reff{hypertweedejac} and \reff{hypereerstejac} by
    induction on $m$.
    For $m=1$, equations \reff{hypertweedejac} and \reff{hypereerstejac}
    reduce to \reff{onehyper} and \reff{twohyper}, respectively.

    In case $m \geq 2$, we use formula \reff{hulpje}, where the
    induction hypothesis enables us to express the right hand polynomial
    $P_{\v n^{(1)}}^{(\v \alpha^{(1)},\beta+n_1)}$ as
    a hypergeometric sum. After exchanging
    the order of summation and differentiation (for $|t|<1$ in
    case of formula (\ref{hypereerstejac})), we apply the formulas
        \begin{eqnarray}
        \label{f1}
        \lefteqn{t^{-\alpha_1}(1-t)^{-\beta}\frac{d^{n_1}}{dt^{n_1}}\left( t^{\alpha_1+n_1+|\v
        k^{(1)}|}(1-t)^{\beta+n_1}\right)} & & \\
        &= & \nonumber
        (\alpha_1+1)_{n_1}\sum_{k_1=0}^{n_1}\frac{(\alpha_1+\beta+n_1+1)_{|\v k|}}
        {(\alpha_1+1)_{|\v k|}}\frac{(\alpha_1+n_1+1)_{|\v k|-k_1}}
        {(\alpha_1+\beta+n_1+1)_{|\v k|-k_1}}
        \frac{(-n_1)_{k_1}t^{|\v k|}}{k_1!}
        \end{eqnarray}
    and
        \begin{equation}
        \label{f2}
        t^{-\alpha_1} \frac{d^{n_1}}{dt^{n_1}}
        t^{\alpha_1+n_1+i}
         =
        \frac{(\alpha_1+1)_{n_1}(\alpha_1+n_1+1)_{i}}{(\alpha_1+1)_{i}}\
        t^{i},
        \end{equation}
    and obtain the right-hand expressions of \reff{hypertweedejac},
    and \reff{hypereerstejac}, respectively.
    It remains to prove the claims \reff{f1} and \reff{f2}, the second one
    being obvious. We observe that the
    left-hand side of \reff{f1} can be transformed using the
    Rodrigues
    formula for $m=1$ and \reff{onehyper}, leading to the
    expression
        \begin{eqnarray*}
        \lefteqn{t^{|\v k^{(1)}|}(\alpha_1+|\v k^{(1)}|+1)_{n_1}
        \f{2}{F}{1}{-n_1,\alpha_1+|\v k^{(1)}|+\beta+n_1+1}{\alpha_1+|\v
        k^{(1)}|+1}{t}} & & \\
        & = &
        (\alpha_1+|\v k^{(1)}|+1)_{n_1}\sum_{k_1=0}^{n_1}\frac{(-n_1)_{k_1}(\alpha_1+|\v k^{(1)}|+\beta+n_1+1)_{k_1}}
        {(\alpha_1+|\v k^{(1)}|+1)_{k_1}}\frac{t^{|\v k|}}{k_1!}\\
        & = &
        (\alpha_1+1)_{n_1}\sum_{k_1=0}^{n_1}\frac{(\alpha_1+\beta+n_1+1)_{|\v k|}}
        {(\alpha_1+1)_{|\v k|}}\frac{(\alpha_1+n_1+1)_{|\v k|-k_1}}
        {(\alpha_1+\beta+n_1+1)_{|\v k|-k_1}}
        \frac{(-n_1)_{k_1}t^{|\v k|}}{k_1!},
        \end{eqnarray*}
    as claimed in \reff{f1}.
\end{proof}
In the above proof we have shown
implicitly that the right hand side of (\ref{hypereerstejac}) is a
polynomial of degree $\leq |\v n|$ in $t$.

With the explicit formula \reff{hypertweedejac} it is possible to
compute the recurrence coefficients by comparing the highest
coefficients in the recurrence relation of order $m+1$, see
\cite{Coussement} for $m=2$. The resulting formulas are quite
complicated, and thus we omit further details.

\subsection{Formal multiple Wilson polynomials}
\label{multwilpolsection}

We now consider $m$ Wilson weights
    \begin{equation}
    \label{wilsonmweights}
    w(\cdot;a,b_j,c,d) ,\quad j=1,\ldots,m, \quad
    \mbox{$b_\ell - b_k \not \in \Z$,\quad  $\ell\neq k$},
    \end{equation}
defined as in \reff{Wilson_weight}, that is, we change only one
parameter (recall the symmetry of the Wilson weights in all four
parameters). Similar as in the scalar case (Wilson polynomials) we
suppose that for $j=1,\ldots,m$
    \begin{equation}
    \label{multiconditions}
    2a,a+b_j,a+c,a+d,2b_j,b_j+c,b_j+d,2c,c+d,2d \notin \{0,-1,-2,\ldots
    \},
    \end{equation}
so that the $m$ Wilson weights have only simple poles, and that
    \begin{equation}
    \label{multicondition2}
    a+b_j+c+d\notin \{0,-1,-2,\ldots \}.
    \end{equation}
As in the scalar case, we write $\mu^{(a,b_j,c,d)}$ for the
resulting measures, where it is possible to choose a joint contour
${\cal C}$ which is the imaginary axis deformed so as to separate
the increasing sequences of poles
($\{a+k\}_{k=0}^\infty,\{b_1+k\}_{k=0}^\infty,\ldots,
\{b_m+k\}_{k=0}^\infty,\{c+k\}_{k=0}^\infty,\{d+k\}_{k=0}^\infty$)
from the decreasing ones
($\{-a-k\}_{k=0}^\infty,\{-b_1-k\}_{k=0}^\infty,\ldots,
\{-b_m-k\}_{k=0}^\infty,\{-c-k\}_{k=0}^\infty,\{-d-k\}_{k=0}^\infty$).

Let us show that the (possibly complex) Wilson measures
$\mu^{(a,b_j,c,d)}$ form a perfect system. The corresponding
multiple orthogonal polynomials will then be referred to as formal
multiple Wilson polynomials. A basic observation in what follows
is that, under some additional conditions, the formal multiple
Wilson polynomials can be written as a Jacobi-Pi\~{n}eiro
transform, similar to \reff{scal_trans}.
    \begin{theorem}
    \label{generaltransform}
    Suppose that \reff{multiconditions} and \reff{multicondition2}
    hold and that $b_\ell - b_k \not \in \Z$, $\ell\neq k$.
    Then the Wilson measures $\mu^{(a,b_j,c,d)}$ for $j=1,\ldots,m$
    form a perfect system.

    If $\Re(a)>0$ and $\Re(c+d)>0$, the formal multiple
    Wilson polynomial with multi-index $\v n$ can be written as
        \begin{equation}
        \label{sim_trans}
         p_{\v n}(t^2; a,\v b,c,d) =
         k_{\v n}\int_{0}^{1}
          P^{(\v \alpha,\beta)}_{\v n}(u)
          w^{(a-1,\beta)}(u)  K(u,t;a,0,c,d)\: du,
        \end{equation}
    for $0<|\Re(t)|<\Re(a)$, where
    $\v \alpha=(a+b_1-1,\ldots,a+b_m-1)=(a-1)  \v e + \v b$ and
    $\beta=c+d-1$.  The normalizing constant
    $k_{\v n}=\v n!\ \Gamma (a+c+|\v n|)\Gamma (a+d+|\v n|)$
    is chosen so that it corresponds with \reff{scal_trans} in the case
    $m=1$ and the kernel $K(u,t;a,b,c,d)$ is defined as in \reff{kern1}.
    \end{theorem}
    Before we prove this theorem we need some technical lemmas.
    \begin{lemma}
    \label{lemmaperfect}
    The system of $m$ measures $\mu_1,\ldots,\mu_m$ is perfect if and only
    if, for every multi-index $\v n$, there exists a polynomial $p_{\v n}$
    of exactly degree $|\v n|$ so that
        \begin{equation}
        \label{equalities}
        \int p_{\v n}(z) z^k \: d\mu_j (z) =0,
        \qquad k=0,\ldots,n_j-1,\quad j=1,\ldots,m,
        \end{equation}
    and
        \begin{equation}
        \label{extratjes}
        \int p_{\v n}(z) z^{n_j} \: d\mu_j (z) \not=0,
        \qquad  j=1,\ldots,m.
        \end{equation}
    In this case, $p_{\v n}$ is the (up to normalization unique)
    multiple orthogonal polynomial of type II with respect to $\v n$.
\end{lemma}
\begin{proof}
    Suppose first that $\mu_1,...,\mu_m$ is perfect, and take as
    $p_{\v n}$ the
    multiple orthogonal polynomial of type II with respect to $\v n$.
    Then it only remains to verify \reff{extratjes}. Suppose the
    contrary, that is, $\int p_{\v n}(z) z^{n_j} \: d\mu_j (z)=0$ for
    some $j$. Then $p_{\v n}$ is also a
    multiple orthogonal polynomial of type II with respect to
    $\v \nu=(n_1,...,n_{j-1},n_j+1,n_{j+1},...,n_m)$, in contrast
    to the normality of the multi-index $\v \nu$.

    We will prove the other implication of this lemma by showing
    by induction on the length $|\v n|$ that $\v n$ is normal.
    The multi-index $\v 0$ is
    always normal, suppose therefore that $\v n$ is of length
    $\geq 1$, with its $j$th component strictly greater than $0$.
    Let $q_{\v n}$ be a multiple orthogonal
    polynomial for $\v n$. If $q_{\v n}$ would have degree strictly less than $|\v
    n|$, then it would be
    a multiple orthogonal polynomial for the multi-index
    $\v \nu=(n_1,...,n_{j-1},n_j-1,n_{j+1},...,n_m)$ with
    length $|\v n|-1$.  By normality of $\v \nu$ we then have that there exists a nonzero
    constant $c$ so that $q_{\v n}=cp_{\v \nu}$.  Thus
    $\int q_{\v n}(z) z^{n_j-1} \: d\mu_j (z)
    \not=0$ by \reff{extratjes}, in contradiction to
    the orthogonality relation \reff{equalities} for $q_{\v n}$.
    As a consequence,
    $q_{\v n}$ has the precise degree $|\v n|$, and thus $\v n$ is
    normal.
\end{proof}
\begin{lemma}
    \label{lemmaperfect2}
    Suppose that the singleton systems $\mu_j$ form a perfect
    system for $j=1,...,m$, with corresponding orthogonal
    polynomials $\{p_{n}^{(j)}\}_n$. Then the system of $m$ measures
    $\mu_1,\ldots,\mu_m$ is perfect if and only
    if, for every multi-index $\v n$, there exists a polynomial $p_{\v n}$
    and scalars $c_{\v n,k}^{(j)}$
    so that
    \begin{equation} \label{eq_lemmaperfect2}
       p_{\v n}(x) = \sum_{\ell=n_j}^{|\v
       n|} c_{\v n,\ell}^{(j)} p_{\ell}^{(j)}(x) , \qquad
       c_{\v n,n_j}^{(j)} \neq 0 , \qquad
       c_{\v n,|\v n|}^{(j)} \neq 0,\qquad j=1,\ldots,m.
    \end{equation}
    In this case, $p_{\v n}$ is the (up to normalization unique)
    multiple orthogonal polynomial of type II with respect to $\v n$.
\end{lemma}
\begin{proof}
    If $\v n$ is normal and $p_{\v n}$ is the corresponding
    multiple orthogonal polynomial, then \reff{eq_lemmaperfect2}
    follows by taking
    \begin{equation} \label{formula_c}
        c_{\v n,\ell}^{(j)} = \frac{\int p_{\v n}(x) p_{\ell}^{(j)}(x)
        d\mu_j(x)}{\int \left(p_{\ell}^{(j)}(x)\right)^2 d\mu_j(x)} ,
    \end{equation}
    where we observe that the denominator is nontrivial according
    to \reff{equalities},\reff{extratjes} for the singleton system
    $\mu_j$. In addition, $c_{\v n,\ell}^{(j)}=0$ for $\ell<n_j$ by
    \reff{equalities}, $c_{\v n,n_j}^{(j)}\neq 0$ by
    \reff{extratjes}, and $c_{\v n,|\v n|}^{(j)}\neq 0$.

    Conversely, \reff{eq_lemmaperfect2} plus the perfectness of
    the singleton system $\mu_j$ implies
    \reff{equalities},\reff{extratjes}, and hence the system $\mu_1,\ldots,\mu_m$ is perfect
    by Lemma \ref{lemmaperfect}.
\end{proof}
    We now prove Theorem \ref{generaltransform} by showing that
    (the analytic extension of) the integral expression
    \reff{sim_trans} is a possible candidate for a formal multiple
    Wilson polynomial.
    \rm \trivlist \item[\hskip \labelsep{\bf Proof of Theorem \ref{generaltransform}. }]
    According to the assumptions \reff{multiconditions} and
    \reff{multicondition2} of Theorem~\ref{generaltransform}, we
    find that
    $\alpha_j,\beta,\alpha_j+\beta+1 \in\C\setminus
    \{-1,-2,\ldots\}$, $j=1,\ldots,m$,
    and $\alpha_\ell - \alpha_k \not \in \Z$ whenever $\ell\neq
    k$.  It follows from Theorem~\ref{perfect_Jacobi_Pinero} that
    the Jacobi-Pi\~{n}eiro system $\mu^{(\alpha_j,\beta)}$,
    $j=1,...,m$, is perfect. From Lemma~\ref{lemmaperfect2} we may
    conclude that there exists scalars $c_{\v n,k}^{(j)}$
    so that
        \[P^{(\v \alpha,\beta)}_{\v n}(t) = \sum_{\ell=n_j}^{|\v n|}
        c_{\v n,\ell}^{(j)} P^{(\alpha_j,\beta)}_{\ell}(t) , \qquad
       c_{\v n,n_j}^{(j)} \neq 0 , \qquad
       c_{\v n,|\v n|}^{(j)} \neq 0, j=1,\ldots,m.\]
    From, e.g., \reff{hypertweedejac} we see that
    Jacobi-Pi\~{n}eiro polynomials are rational functions in each of
    the parameters $\alpha_j$ or $\beta$. Taking into account
    \reff{beta} and \reff{formula_c}, we may conclude that any of
    the coefficients $c_{\v n,\ell}^{(j)}$ is a meromorphic function
    in each of the parameters $\alpha_j$ or $\beta$.
    We now introduce
        \begin{eqnarray*}
         q_{\v n}(t^2; a,\v b,c,d) & = & k_{\v n}
         \int_{0}^{1} P^{(\v \alpha,\beta)}_{\v n}(u)
          w^{(a-1,\beta)}(u) K(u,t;a,0,c,d)\: du.
        \end{eqnarray*}
    This function is well defined for $0<|\Re(t)|<\Re(a)$ if $\Re(a)>0$
    and $\Re(c+d)>0$.
    However, using \reff{scal_trans}, we
    obtain for every $j=1,\ldots,m$ that
        \begin{eqnarray*}
        q_{\v n}(t^2; a,\v b,c,d) & = &
        k_{\v n}\sum_{\ell=n_j}^{|\v n|}
        c_{\v n,\ell}^{(j)}
        \int_{0}^{1}P^{(\alpha_j,\beta)}_{\ell}(x) w^{(\alpha_j,\beta)}(u) K(u,t;a,b_j,c,d)\: du,
        \\ & = &
        \sum_{\ell=n_j}^{|\v n|}
        \frac{\v n!}{\ell !}
        \frac{\Gamma(a+c+|\v n|)\Gamma(a+d+|\v n|)}{\Gamma(a+c+\ell)\Gamma(a+d+\ell)}
        \ c_{\v n,\ell}^{(j)} \
        p_{\ell}(t^2; a,b_j,c,d) ,
        \end{eqnarray*}
    and thus
        \begin{eqnarray}
            \nonumber
            q_{\v n}(t^2; a,\v b,c,d) & = &
            \sum_{\ell=n_j}^{|\v n|} \frac{\v n!}{\ell!}
            (a+c+\ell)_{|\v n|-\ell}\:(a+d+\ell)_{|\v n|-\ell}\
            c_{\v n,\ell}^{(j)} \ p_{\ell}(t^2; a,b_j,c,d)\\
            & = & \label{def_multiple_Wilson}
            \sum_{\ell=n_j}^{|\v n|}
            d_{\v n,\ell}^{(j)} \ p_{\ell}(t^2; a,b_j,c,d).
        \end{eqnarray}
    Observing that the expressions on the right-hand side
    of \reff{def_multiple_Wilson}
    are polynomials in $t$ and meromorphic in any of the
    parameters
    $a,\v b,c,d$, we see that the right hand expression
    of \reff{def_multiple_Wilson} is well defined
    and independent of $j$
    for any choice of $t$ and of the parameters $a,\v b,c,d$, as long as
    \reff{multiconditions} and \reff{multicondition2}
    hold and $b_\ell - b_k \not \in \Z$, $\ell\neq k$. Moreover, the new coefficients
    $d_{\v n,n_j}^{(j)}$ and for $d_{\v n,|\v n|}^{(j)}$ are different from zero.

    Thus, $q_{\v n}(t^2; a,\v b,c,d)$ defined by \reff{def_multiple_Wilson}
    is a suitable candidate for
    a formal multiple Wilson polynomial, and the system of Wilson
    measures is perfect by Lemma~\ref{lemmaperfect2}.
{\hspace*{\fill}$\Box$\endtrivlist}
We now want to deduce explicit expressions for the formal multiple
Wilson polynomials based on the explicit expressions
\reff{hypertweedejac} and \reff{hypereerstejac} for the
Jacobi-Pi\~{n}eiro polynomials. Here we use the Kamp\'e de
F\'eriet series \cite{Srivastava}
    \begin{equation}
    \label{Feriet}
    \g{F}{q:q_1;q_2}{p:p_1;p_2}{\v f:\v g;\v h}{\vec{\phi}:\vec{\psi};\vec{\xi}}{x,y}
    := \sum_{k=0}^\infty
    \frac{\prod\limits_{\ell=1}^{p}(f_\ell)_{k}}
         {\prod\limits_{\ell=1}^{q}(\phi_\ell)_{k}}
    \sum_{j=0}^k
    \frac{\prod\limits_{\ell=1}^{p_1}(g_\ell)_{k-j}
          \prod\limits_{\ell=1}^{p_2}(h_\ell)_{j}}
         {\prod\limits_{\ell=1}^{q_1}(\psi_\ell)_{k-j}
          \prod\limits_{\ell=1}^{q_2}(\xi_\ell)_{j}}
    \frac{x^{k-j}}{(k-j)!}\frac{y^j}{j!},
    \end{equation}
which are a generalization of the 4 Appell series in 2 variables.
Notice that, for $p=q=0$, the Kamp\'e de F\'eriet series is a
product of two hypergeometric series. Also, in the case $m=2$, our
functions $\mathcal{M}_{q,\v n}^{p;r}$ defined in
(\ref{nieuweformule}) are (finite) Kamp\'e de F\'eriet series
\begin{equation} \label{link_M_Kampe}
    \g{\mathcal{M}}{q,(n_1,n_2)}{p;r}{\v f ;g_1:\cdots :g_r}{\v \phi ; \psi_1:\cdots :\psi_r}{
    (x_1,x_2)} =
    \g{F}{q:0;r}{p:1;r+1}{\v
    f:(-n_1);(-n_2,g_1,...,g_r)}{\vec{\phi}:();(\psi_1,...,\psi_r)}{x_1,x_2}.
\end{equation}
In what follows the parameters in the Kamp\'e de F\'eriet series
will always be chosen so that the sum in (\ref{Feriet}) is finite,
and hence we are not concerned with convergence problems.
\begin{corollary}\label{representation_wilson}
    Let $\v e=(1,\ldots,1)$ be a multi-index of length $m$,
    $s(\v n)=(n_1,n_1+n_2,\ldots,|\v n|)$ and $\sigma_j=a+b_j+c+d-1$, $j=1,\ldots,m$.
    Denote by $\v v^{(j)}$ the vector $\v v$
    without the $j$th component.  For the
    multiple Wilson polynomials we have the hypergeometric
    representations
        \begin{eqnarray}
        \label{hypereerstewil}
        & & p_{\v n}(t^2; a,\v b,c,d)=
        (a\v e + \v b)_{\v n}(a+c)_{|\v n|}(a+d)_{|\v n|}\\
        \nonumber
        & & \qquad \qquad\times \g{\mathcal{M}}{3,\v n}{3;2}{(a-t,a+t,\sigma_1+n_1);(a\v e+ \v b+\v n)^{(m)}:
        (\v \sigma+s(\v n))^{(1)}}
        {(a+c,a+d,a+b_1) ;(a\v e+\v b)^{(1)}:(\v \sigma+s(\v n))^{(m)}}{\v e}
        \end{eqnarray}
    and
        \begin{eqnarray}
        \label{hypertweedewil}
        & & p_{\v n}(t^2; a,\v b,c,d) =
            (a\v e + \v b)_{\v n}(a+c)_{|\v n|}(a+d)_{|\v n|}\\
        \nonumber
        & & \qquad \qquad \times
            \g{F}{2:0;m}{2:1;m+1}{(a-t,a+t):(c+d-1);(a \v e+ \v b +\v n ,1-c-d-|\v n|)}
                                {(a+c,a+d):();(a \v e+\v b)}{1,1}
        \end{eqnarray}
    where $(a\v e+\v b)_{\v n} = \prod_{j=1}^m(a+b_j)_{n_j}$.
\end{corollary}
\begin{proof}
    First note that $(1-x)^{-\beta}=\sum_{\ell=0}^\infty (\beta)_l\frac{x^\ell}{\ell!}$,
    which converges in the unit disk, so that the expression
    \reff{hypereerstejac} can then be written as
         \begin{equation}
        \label{uitgewerkt}
        P_{\v n}^{(\v \alpha,\beta)}(t)=
        \frac{(\v \alpha+\v e)_{\v n}}{\v n!}
        \g{F}{0:0;m}{0:1;m+1}{-:(\beta);(\v \alpha+\v n + \v e,-\beta-|\v
        n|)}{-:-;(\v \alpha + \v e)}{t,t}.
        \end{equation}
    Then start from the Jacobi-Pi\~{n}eiro transform \reff{sim_trans}
    and replace the  Jacobi-Pi\~{n}eiro polynomial $P_{\v n}^{(\v
    \alpha,\beta)}$ by its explicit expressions \reff{hypertweedejac} and \reff{uitgewerkt},
    respectively.  Since the sums are finite,
    we can interchange the integral with the sums.
    Applying \reff{bewijstransform} then completes the proof.
\end{proof}
For $m=2$, we may apply (\ref{link_M_Kampe}) to
(\ref{hypereerstewil}), leading to a representation as a Kamp\'e
de F\'eriet series of type $F_{3:0;2}^{3:1;3}$. It seems to be
non-trivial to derive from this formula the representation
(\ref{hypertweedewil}) as a Kamp\'e de F\'eriet series of type
$F_{2:0;2}^{2:1;3}$.

\section{Limit relations}
\label{multAskeysection}

In this section we consider some cases in which the orthogonality
conditions of the formal multiple Wilson polynomials reduce to
orthogonality conditions with respect to a positive measure on the
real line. We then recover the multiple Wilson and multiple Racah
polynomials. Next we use the limit relations between the
orthogonal polynomials in the Askey table \cite{Koekoek} to obtain
some new examples of multiple orthogonal polynomials and some
known examples. In particular we look at what happens with the
explicit expressions \reff{hypereerstewil} and
\reff{hypertweedewil} applying these limit relations.  Most of
these examples are known to be AT systems (which implies that
every multi-index is normal).

\subsection{Multiple Wilson}
With some restrictions on the parameters the orthogonality
conditions of the formal multiple Wilson polynomials reduce to
real orthogonality conditions with respect to positive measures on
the real line. Let $b_j>0$, $j=1,\ldots,m$, $b_\ell - b_k \not \in
\Z$ whenever $\ell\neq k$, $\Re(a),\Re(c),\Re(d)>0$ and $a,c,d$ be
real except for a conjugate pair. In that case the imaginary axis
can be taken as the contour $\cal C$. The multiple Wilson
polynomials
    \begin{equation}
    \label{realmultWilson}
    W_{\v n}(x^2;a,\v b,c,d):=p_{\v n}(-x^2;a,\v b,c,d)
    \end{equation}
then satisfy the real orthogonality relations
    \begin{equation*}
    \int_0^{\infty} (x^2)^k \: W_{\v n}(x^2;a,\v b,c,d)
    \left|\frac{\Gamma(a+ix)\Gamma(b_j+ix)\Gamma(c+ix)\Gamma(d+ix)}{\Gamma(2ix)}\right|^2
    dx = 0,
    \end{equation*}
$k=0,\ldots,n_j-1$, $j=1,\ldots,m$.  If $a<0$, $a+b_j>0$,
$j=1,\ldots,m$, and $a+c,a+d$ are positive or a pair of complex
conjugates with positive real parts, then we obtain the same
orthogonality conditions but with some extra positive point
masses.

\subsection{Multiple Racah}
As in the scalar case it is also possible to obtain a purely
discrete orthogonality. The multiple Racah polynomials $R_{\v
n}(\cdot;\v \alpha, \beta,\gamma,\delta)$ (where we only change
the parameter $\alpha$ with $\alpha_\ell - \alpha_k \not \in \Z$
whenever $\ell\neq k$) satisfy the discrete orthogonality
    \[
    \sum_{x=0}^N\frac{(\alpha_j+1)_x(\gamma+1)_x(\beta+\delta+1)_x(\gamma+\delta+1)_x((\gamma+\delta+3)/2)_x}
    {(-\alpha_j+\gamma+\delta+1)_x(-\beta+\gamma+1)_x((\gamma+\delta+1)/2)_x(\delta+1)_xx!}\
    R_{\v n}(\lambda(x);\v \alpha,\beta,\gamma,\delta)\:(\lambda(x))^k=0,
    \]
$k=0,\ldots,n-1$, $j=1,\ldots,m$, where
    \[\lambda(x)=x(x+\gamma+\delta+1) \qquad \mbox{and} \qquad \beta+\delta+1=-N \mbox{\ \   or\ \ } \gamma+1=-N.\]
They can be found from the polynomials
    $p_{\v n}(t^2;a,\v b,c,d)/( {(a\v e + \v b)_{\v n}(a+c)_{|\v n|}(a+d)_{|\v n|}})$
by the substitution $t\to x+a$ and the change of variables
$\alpha_j=a+b_j-1, \beta=c+d-1, \gamma=a+d-1, \delta=a-d$. For the
multiple Racah polynomials we then have the expressions
    \begin{eqnarray}
    \nonumber
    & & R_{\v n}(\lambda(x);\v \alpha,\beta,\gamma,\delta)\\
    \nonumber & & =
    \g{\mathcal{M}}{3,\v n}{3;2}{(-x,x+\gamma+\delta+1,\alpha_1+\beta+n_1+1);(\v \alpha+\v n+\v e)^{(m)}:
    (\v \alpha+s(\v n)+(\beta+1)\v e)^{(1)}}
    {(\beta+\delta+1,\gamma+1,\alpha_1+1) ;(\v \alpha+\v e)^{(1)}:
    (\v \alpha+s(\v n)+(\beta+1)\v e)^{(m)}}{\v e}\\
    \nonumber
    & & =  \g{F}{2:0;m}{2:1;m+1}{(-x,x+\gamma+\delta+1):(\beta);(\v \alpha+ \v n +\v e,-\beta-|\v n|)}
    {(\beta+\delta+1,\gamma+1):();(\v \alpha +\v e)}{1,1}.
    \end{eqnarray}
In this case, necessary and sufficient conditions on the
parameters to have positive weights are quite messy.
    \begin{remark}
    Remember that the Wilson weight is symmetric in the four
    parameters so that we can switch these parameters in the
    change of variables.  We then obtain multiple Racah polynomials (of type II)
    where we change other parameters.  However, all these cases are related to each other.
    For example we have multiple Racah polynomials $R_{\v n}(\cdot;\alpha,\v \beta,\gamma,\delta)$ where we only
    change the parameter $\beta$ in the weights ($\beta_\ell - \beta_k \not \in \Z$ whenever $\ell\neq k$).  In that case we
    have that $\alpha+1=-N$ or $\gamma+1=-N$.
    As a second example it is possible to change the parameters $\beta,\gamma$ and $\delta$
    in such a way that $\delta_j+\gamma_j$ and $\delta_j+\beta_j$
    don't change and $\gamma_\ell - \gamma_k \not \in \Z$ whenever $\ell\neq k$.
    Here we assume that $\alpha+1=-N$ or $\beta_j+\delta_j+1=-N$.
    We then denote these multiple Racah polynomials by
    $R_n(\cdot;\alpha,\v \beta,\v \gamma,\v \delta)$.  However, we don't find another
    family of polynomials because
        \begin{equation}
        \label{tussenRacah1}
        R_{\v n}(\lambda(x);\alpha,\v \beta,\gamma,\delta)=
        R_{\v n}(\lambda(x);\v \beta+\delta\v e,\alpha-\delta,\gamma,\delta)
        \end{equation}
    and
        \begin{equation}
        \label{tussenRacah}
        R_{\v n}(\lambda(x);\alpha,\v \beta,\v \gamma,\v \delta)=
        R_{\v n}(\lambda(x);\v
        \gamma,\alpha+\beta_j-\gamma_j,\alpha,\gamma_j+\delta_j-\alpha).
        \end{equation}
    These relations will help us in some of the examples of the
    subsections below to find explicit expressions for the
    polynomials.
    \end{remark}

\subsection{Some new examples}

\subsubsection{Multiple continuous dual Hahn}
Let $b_\ell - b_k \not \in \Z$ whenever $\ell\neq k$.  The
multiple continuous dual Hahn polynomials satisfy the
orthogonality conditions of the multiple Wilson polynomials where
we let $d\to +\infty$ (after dividing by $\Gamma(d)^2$). Similarly
we obtain real orthogonality conditions with respect to a positive
measure if $b_j>0$ and $a,c$ are positive or a pair of complex
conjugates with positive real parts.  We then denote the $\v n$th
multiple continuous dual Hahn polynomial by $S_{\v n}(\cdot ;a,\v
b,c)$. These polynomials satisfy the orthogonality conditions
    \begin{equation*}
    \int_0^{\infty} (x^2)^k \: S_{\v n}(x^2;a,\v b,c)
    \left|\frac{\Gamma(a+ix)\Gamma(b_j+ix)\Gamma(c+ix)}{\Gamma(2ix)}\right|^2
    dx = 0,
    \end{equation*}
$k=0,\ldots,n_j-1$, $j=1,\ldots,m$.  It is clear that
    \begin{equation}
    S_{\v n}(x^2;a,\v b,c)=\lim_{d\to +\infty}\frac{W_{\v n}(x^2;a,\v
    b,c,d)}{(a+d)_{|\v n|}}
    \end{equation}
so that the multiple continuous dual Hahn polynomials have the
explicit expressions
    \begin{eqnarray*}
        p_{\v n}(x^2; a,\v b,c) & = &
        (a\v e + \v b)_{\v n}(a+c)_{|\v n|}
        \g{\mathcal{M}}{2,\v n}{2;1}{(a-ix,a+ix);(a\v e+ \v b+\v n)^{(m)}}
        {(a+c,a+b_1) ;(a\v e+\v b)^{(1)}}{\v e},\\
        & = &
        (a\v e + \v b)_{\v n}(a+c)_{|\v n|}
            \g{F}{1:0;m}{2:0;m}{(a-ix,a+ix):();(a \v e+ \v b +\v n )}
                                {(a+c):();(a \v e+\v b)}{1,-1}.
    \end{eqnarray*}

\subsubsection{Multiple dual Hahn}

Consider $\gamma_j,\delta_j$, $j=1,\ldots,m$, so that
$\gamma_j,\delta_j>-1$ or $\gamma_j,\delta_j<-N$ for each $j$ and
that $\gamma_j+\delta_j$ is independent of $j$.  Suppose also that
$\gamma_\ell - \gamma_k \not \in \Z$ whenever $\ell\neq k$. The
multiple dual Hahn polynomials, denoted by $R_{\v n}(\cdot; \v
\gamma,\v \delta,N)$, satisfy the system of discrete orthogonality
conditions
    \[
    \sum_{x=0}^N\frac{(2x+\gamma_j+\delta_j+1)(\gamma_j+1)_x(-N)_xN!}
    {(-1)^x(x+\gamma_j+\delta_j+1)_{N+1}(\delta_j+1)_xx!}\
    R_{\v n}(\lambda(x);\v \gamma,\v \delta,N)\:(\lambda(x))^k=0,
    \]
$k=0,\ldots,n-1$, $j=1,\ldots,m$, where $\lambda(x)=
x(x+\gamma+\delta+1)$.  The multiple dual Hahn polynomials are
related to the multiple Racah polynomials like
    \begin{eqnarray}
    \nonumber R_{\v n}(\lambda(x);\v \gamma,\v \delta,N)
    & = & \lim_{\alpha \to +\infty} R_{\v n}(\lambda(x);
    \alpha,-\v \delta-(N+1)\v e,\v \gamma,\v \delta)\\
    & = & \lim_{\alpha \to +\infty} R_{\v n}(\lambda(x);
    \v \gamma,\alpha-\gamma_j-\delta_j-N-1,
    \alpha,\gamma_j+\delta_j-\alpha),
    \end{eqnarray}
where we use \reff{tussenRacah}.  They then have the explicit
expressions
    \begin{eqnarray*}
    R_{\v n}(\lambda(x);\v \gamma,\v \delta,N)
    & = &
    \g{\mathcal{M}}{2,\v n}{2;1}{(-x,x+\gamma_j+\delta_j+1);(\v \gamma+\v n+\v e)^{(m)}}
    {(-N,\gamma_1+1) ;(\v \gamma+\v e)^{(1)}}{\v e}\\
    & = &  \g{F}{1:0;m}{2:0;m}{(-x,x+\gamma_j+\delta_j+1):();(\v \gamma+ \v n +\v e)}
    {(-N):();(\v \gamma +\v e)}{1,-1}.
    \end{eqnarray*}

\subsubsection{Multiple Meixner-Pollaczek}

The multiple Meixner-Pollaczek polynomials $P_{\v
n}^{(\lambda)}(\cdot;\v \phi)$ are multiple orthogonal polynomials
(of type II) associated with the system of weights
$e^{(2\phi_j-\pi)x}|\Gamma(\lambda+ix)|^2$, $\lambda>0$,
$0<\phi_j<\pi$, $j=1,\ldots,m$, where the $\phi_1,\ldots,\phi_m$
are different.  These weights form an AT system on the positive
real axis \cite[p.141]{Nikishin}.  The multiple Meixner-Pollaczek
polynomials satisfy the conditions
    \[
    \int_0^{\infty} x^k \: P_{\v n}^{(\lambda)}(x;\v \phi)
    e^{(2\phi_j-\pi)x}|\Gamma(\lambda+ix)|^2\:
    dx = 0,
    \]
$k=0,\ldots,n_j-1$, $j=1,\ldots,m$.  Similar as
\cite[(2.3.1)]{Koekoek} it is easy to check that
    \begin{equation}
    P_{\v n}^{(\lambda)}(x;\v \phi)=\lim_{t\to +\infty}\frac{S_{\v n}((x-t)^2;\lambda+it,t \cot\v\phi,\lambda-it)}
    {(t\csc \v \phi)_{\v n}\: \v n!},
    \end{equation}
where $t \cot\v\phi=(t\cot\phi_1,\ldots,t\cot\phi_m)$ and $t\csc
\v \phi=(t\csc \phi_1,\ldots,t\csc \phi_m)$.  The multiple
Meixner-Pollaczek polynomials then have the explicit expression
    \begin{equation*}
    P_{\v n}^{(\lambda)}(x;\v \phi)  =
    \frac{(2\lambda)_{|\v n|}\prod_{j=1}^m e^{in_j\phi_j}}{\v n!}
    \g{\mathcal{M}}{1,\v n}{1;0}{(\lambda+ix);-}
    {(2\lambda) ;-}{\v e-e^{-2i\v \phi}},
    \end{equation*}
where $e^{-2i\v \phi}=(e^{-2i\phi_1},\ldots,e^{-2i\phi_m})$. Here
we do not have a Kamp\'e de F\'eri\'et representation such as in
\reff{hypertweedewil}.

\subsubsection{Formal multiple continuous Hahn}
Similar as in \cite[(2.1.2)]{Koekoek} we can use the limit
relation
    \begin{equation}
    \label{limitcontHahn}
    P_{\v n}(t;a,\v b,c,d)=\lim_{y\to \infty}
    \frac{p_{\v n}((t+y)^2; a-y,\v b+y\v e,c-y,d+y)}{(a+c-2y)_{|\v n|}\ \v n!}
    \end{equation}
in order to find the formal continuous Hahn polynomials.  They
have the explicit expressions
        \begin{eqnarray*}
        P_{\v n}(t; a,\v b,c,d) & = &
        (a\v e + \v b)_{\v n}(a+d)_{|\v n|}
        \g{\mathcal{M}}{2,\v n}{2;2}{(a+t,\sigma_1+n_1);(a\v e+ \v b+\v n)^{(m)}:
        (\v \sigma+s(\v n))^{(1)}}
        {(a+d,a+b_1) ;(a\v e+\v b)^{(1)}:(\v \sigma+s(\v n))^{(m)}}{\v   e}\\
        & = & (a\v e + \v b)_{\v n}(a+d)_{|\v n|}\\
        & & \qquad \quad \times
            \g{F}{1:0;m}{1:1;m+1}{(a+t):(c+d-1);(a \v e+ \v b +\v n ,1-c-d-|\v n|)}
                                {(a+d):();(a \v e+\v b)}{1,1},
        \end{eqnarray*}
where $\sigma_j=a+b_j+c+d-1$, $j=1,\ldots,m$. If the parameters
satisfy \reff{multiconditions} and \reff{multicondition2} and
$b_\ell - b_k \not \in \Z$ whenever $\ell\neq k$, then these
polynomials satisfy the orthogonality conditions
        \begin{equation}
        \label{conthahnorthogmulti}
        \int_{\cal C}P_{\v n}(t; a,\v b,c,d) \Gamma(a+t)\Gamma(b_j-t)\Gamma(c+t)\Gamma(d-t) t^k \: dt
         = 0 ,
         \end{equation}
$k=0,1,\ldots,n_j-1$, $j=1,\ldots,m$, where ${\cal C}$ is a
contour which is the imaginary axis deformed so as to separate the
increasing sequences of poles ($\{b_1+k\}_{k=0}^\infty,\ldots,
\{b_m+k\}_{k=0}^\infty,\{d+k\}_{k=0}^\infty$) from the decreasing
ones ($\{-a-k\}_{k=0}^\infty,\{-c-k\}_{k=0}^\infty$).  In the
scalar case ($m=1$) it is possible to obtain real orthogonality
relations with respect to a positive measure if we suppose
$\Re(a),\Re(b),\Re(c),\Re(d)>0$ and $a= \bar{b}$, $c= \bar{d}$.
This is not possible in the multiple case.  Therefore it would be
nice if we could find another family of multiple Wilson
polynomials in the sense that we change two parameters forming a
pair of complex conjugates.

\subsection{Some classical discrete multiple orthogonal polynomials}
In this section we obtain hypergeometric formulas for the
classical discrete examples of multiple orthogonal polynomials of
type II, introduced in \cite{Arvesu}, which are all examples of AT
systems. In particular we use the limit relations between these
polynomials and the Racah polynomials \cite{Koekoek}.  Their
$\mathcal{M}_{q,\v n}^{p;r}$ representation is already known in
the cases $m=1,2$. The explicit expression in terms of a Kamp\'e
de F\'eri\'et series is new (if it exists).

\begin{itemize}
\item {\bf Multiple Hahn:} These multiple orthogonal polynomials
(of type II) satisfy orthogonality conditions with respect to $m$
hypergeometric distributions
    \[\mu_j=\sum_{k=0}^N\frac{(\alpha_j+1)_k}{k!}\frac{(\beta+1)_{N-k}}{(N-k)!}\delta_k,\quad \alpha_j>-1,
    \ \beta>-1,\]
$\alpha_l-\alpha_k\notin \{0,1,\ldots,N-1\}$, $\ell\not=k$, on the
integers $0,\ldots,N$.  They can be found from the multiple Racah
polynomials taking $\gamma+1=-N$ and $\delta\to +\infty$, so that
    \begin{eqnarray*}
    Q_{\v n}^{\v \alpha;\beta;N}(x)
    & = &
    \g{\mathcal{M}}{2,\v n}{2;2}{(-x,\alpha_1+\beta+n_1+1);(\v \alpha+\v n+\v e)^{(m)}:
    (\v \alpha+s(\v n)+(\beta+1)\v e)^{(1)}}
    {(-N,\alpha_1+1) ;(\v \alpha+\v e)^{(1)}:
    (\v \alpha+s(\v n)+(\beta+1)\v e)^{(m)}}{\v e}\\
    & = &   \g{F}{1:0;m}{1:1;m+1}{(-x):(\beta);(\v \alpha+ \v n +\v e,-\beta-|\v n|)}
    {(-N):();(\v \alpha +\v e)}{1,1}.
    \end{eqnarray*}
Changing only the parameter $\beta$ does not give another family
of polynomials because of $Q_{\v n}^{\alpha;\v \beta;N}(x)=C\:
Q_{\v n}^{\v \beta;\alpha;N}(N-x)$ with $C$ some constant
(depending on $\v n,\alpha$ and $\v \beta$).  However, we will
need an explicit formula in powers of $x$ for these polynomials to
reach multiple Meixner I and multiple Laguerre II. Using
\reff{tussenRacah1} and the limits we mentioned above, we find
that
    \begin{eqnarray*}
    Q_{\v n}^{\alpha;\v \beta;N}(x)
    & = &
    \g{\mathcal{M}}{2,\v n}{2;2}{(-x,\alpha+\beta_1+n_1+1);
    (\v \beta+s(\v n)+(\alpha+1)\v e)^{(1)}}
    {(-N,\alpha+1) ;
    (\v \beta+s(\v n)+(\alpha+1)\v e)^{(m)}}{\v e}.
    \end{eqnarray*}
Applying these limits to the Kamp\'e de F\'eri\'et representation
does not work.
\item {\bf Multiple Meixner I:} In this case we consider $m$
negative binomial distributions
    \[\mu_j=\sum_{k=0}^\infty \frac{(\beta)_k\: c_j^k}{k!}\delta_k,\qquad 0<c_j<1,\ \beta>0,\]
with all the $c_j$, $j=1\ldots,m$, different. We get these
polynomials from the multiple Hahn polynomials $Q_{\v
n}^{\alpha;\v \beta;N}$ replacing $\alpha=\beta-1$,
$\beta_j=N\frac{1-c_j}{c_j}$ and letting $N\to +\infty$. We then
obtain
    \begin{eqnarray*}
    M_{\v n}^{\beta;\v c}(x)
    & = &
    \g{\mathcal{M}}{1,\v n}{1;0}{(-x);()}
    {(\beta) ;()}{\frac{\v c-\v e}{\v c}},
    \end{eqnarray*}
where $\frac{\v c-\v e}{\v
c}=\left(\frac{c_1-1}{c_1},\ldots,\frac{c_m-1}{c_m}\right)$.
\item {\bf Multiple Meixner II:} In the case of multiple Meixner
II polynomials we only change the parameter $\beta$ in the
negative binomial distributions, so that
    \[\mu_j=\sum_{k=0}^\infty \frac{(\beta_j)_k\: c^k}{k!}\delta_k,\qquad 0<c<1,\ \beta_j>0,\]
with $\beta_\ell - \beta_k \not \in \Z$ whenever $\ell\neq k$.
Taking $\alpha_j=\beta_j-1$, $\beta=N\frac{1-c}{c}$ and letting
$N\to +\infty$ in the explicit formulas for the multiple Hahn
polynomials $Q_{\v n}^{\v \alpha;\beta;N}$, we obtain
    \begin{eqnarray*}
    M_{\v n}^{\v \beta;c}(x)
    & = &
    \g{\mathcal{M}}{1,\v n}{1;1}{(-x);(\v \beta+\v n)^{(m)}}
    {(\beta_1) ;(\v \beta)^{(1)}}{\frac{c-1}{c}\v e}\\
    & = &   \g{F}{0:0;m}{1:0;m}{(-x):();(\v \beta+ \v n)}
    {():();(\v \beta)}{\frac{c-1}{c},\frac{1-c}{c}}.
    \end{eqnarray*}
\item {\bf Multiple Kravchuk:}  These polynomials satisfy the orthogonality conditions \reff{system}
with the $m$ binomial distributions
    \[\mu_j=\sum_{k=0}^N{N\choose k}p_j^k(1-p_j)^{N-k}\delta_k,\qquad 0<p_j<1,\]
where all the $p_j$, $j=1\ldots,m$, are different. These
polynomials are related to the multiple Hahn polynomials $Q_{\v
n}^{\v \alpha;\beta;N}$ replacing $\beta= t$, $\alpha_j\to
\frac{p_j}{1-p_j}t$ and letting $t\to +\infty$.  We then get
    \begin{eqnarray*}
    K_{\v n}^{\v p;N}(x)
    & = &
    \g{\mathcal{M}}{1,\v n}{1;0}{(-x);()}
    {(-N) ;()}{\frac{1}{\v p}},
    \end{eqnarray*}
where $\frac{1}{\v p
}=\left(\frac{1}{p_1},\ldots,\frac{1}{p_m}\right)$.
\item {\bf Multiple Charlier:}
In the case of multiple Charlier we consider $m$ Poisson
distributions
    \[\mu_j=\sum_{k=0}^\infty \frac{a_j^k}{k!}\delta_k,\qquad a_j>0,\]
with all the $a_j$, $j=1\ldots,m$, different.   The corresponding
multiple orthogonal polynomials (of type II) can be found from the
multiple Meixner I polynomials taking $c_j=\frac{a_j}{a_j+\beta}$
and letting $\beta\to +\infty$. The multiple Charlier polynomials
then have the explicit expression
    \begin{eqnarray*}
    C_{\v n}^{\v a}(x)
    & = &
    \g{\mathcal{M}}{0,\v n}{1;0}{(-x);()}
    {() ;()}{-\frac{1}{\v a}},
    \end{eqnarray*}
where $\frac{1}{\v
a}=\left(\frac{1}{a_1},\ldots,\frac{1}{a_m}\right)$.
\end{itemize}

\subsection{Some classical continuous multiple orthogonal polynomials}
In this section we recall some classical continuous examples of
multiple orthogonal polynomials of type II where the measures (or
weight functions) form an AT system and obtain hypergeometric
formulas for these polynomials.  Their $\mathcal{M}_{q,\v
n}^{p;r}$ representation is already known in the cases $m=1,2$.
The explicit expression in terms of a Kamp\'e de F\'eri\'et series
is new (if it exists). For an overview of these polynomials and
their properties we recommend \cite{Coussement,Apt}.

\begin{itemize}
\item {\bf Jacobi-Pi\~{n}eiro:}
In Subsection \ref{defjacpin} we recalled the Jacobi-Pi\~{n}eiro
polynomials $P_{\v n}^{\v \alpha,\beta}$, which, in the case
$\alpha_j,\beta>-1$, are the multiple orthogonal polynomials (of
type II) with respect to the Jacobi weights
$w^{\alpha_j,\beta}(x)=x^{\alpha_j}(1-x)^{\beta}$, $j=1,\ldots,m$,
on the interval $[0,1]$.  Here $\alpha_\ell - \alpha_k \not \in
\Z$ whenever $\ell\neq k$.  Similar as in the multiple Hahn case
we have that $P_{\v n}^{(\alpha,\v \beta)}(x)=(-1)^{|\v n|} P_{\v
n}^{(\v \beta,\alpha)}(1-x)$.  So, changing only the parameter
$\beta$ does not give another family of polynomials.  For these
polynomials we have
        \begin{eqnarray*}
        P_{\v n}^{(\alpha,\v \beta)}(x) & = &
        \lim_{N\to +\infty}\frac{(\alpha+1)_{|\v n|}}{\v n!}Q_{\v n}^{\alpha;\v \beta:N}(N x)\\
        & = &
        \frac{(\alpha+1)_{|\v n|}}{\v n!}
        \g{\mathcal{M}}{1,\v n}{1;1}{(\alpha+\beta_1+n_1+1);
        (\v \beta+s(\v n)+ (\alpha+1)\v e)^{(1)}}
        {(\alpha+1) ;(\v \beta+s(\v n)+(\alpha+1)\v e)^{(m)}}{x\v
        e}.
        \end{eqnarray*}
\item {\bf Multiple Laguerre I:}
The multiple Laguerre I polynomials $L_{\v n}^{\v \alpha}$ are
orthogonal on $[0,\infty)$ with respect to the $m$ weights
$w_j(x)=x^\alpha_je^{-x}$, where $\alpha_j>-1$, $j=1,\ldots,m$,
and $\alpha_\ell - \alpha_k \not \in \Z$ whenever $\ell\neq k$.
They can be found from the Jacobi-Pi\~{n}eiro polynomials $P_{\v
n}^{\v \alpha,\beta}$ substituting $x\to \frac{x}{\beta}$ and
letting $\beta \to \infty$.  We then obtain the hypergeometric
expressions
        \begin{eqnarray*}
        L_{\v n}^{\v \alpha}(x) & = &
        \frac{(\v \alpha+\v e)_{\v n}}{\v n!}
        \g{\mathcal{M}}{1,\v n}{0;1}{();(\v \alpha+\v n+\v e)^{(m)}}
        {(\alpha_1+1) ;(\v \alpha+\v e)^{(1)}}{x\v
        e}\\
         & = &
        \frac{(\v \alpha+\v e)_{\v n}}{\v n!} \ e^x
        \f{m}{F}{m}{\v \alpha+\v n + \v e}{\v \alpha + \v e}{-x}.
        \end{eqnarray*}
\item {\bf Multiple Laguerre II:}
In this case the polynomials $L_{\v n}^{(\alpha,\v c)}$ have the
orthogonality conditions \reff{system} with respect to the weight
functions $w_j(x)=x^\alpha e^{-c_jx}$, where $\alpha>-1$, $c_j>0$,
$j=1,\ldots,m$, and all the $c_j$ different.  They can be obtained
from the Jacobi-Pi\~{n}eiro polynomials $P_{\v n}^{(\alpha,\v
\beta)}$ by the substitutions $x\to \frac{x}{t}$, taking
$\beta_j=c_jt$ and letting $t\to \infty$.  We then get
        \begin{eqnarray*}
        L_{\v n}^{(\alpha,\v c)}(x) & = &
        \frac{(\alpha+1)_{|\v n|}}{\v n!}
        \g{\mathcal{M}}{1,\v n}{0;0}{();
        ()}{(\alpha+1) ;()}{x\v c}.
        \end{eqnarray*}
\item {\bf Multiple Hermite:}
In the multiple Hermite case we consider the multiple orthogonal
polynomials (of type II) $H_{\v n}^{\v c}$ with respect to the
weights $w_j(x)=e^{-x^2+c_jx}$, $j=1,\ldots,m$, on
$(-\infty,+\infty)$. Here the $c_j$ are different real numbers.
These polynomials can be found from the Jacobi-Pi\~{n}eiro
polynomials $P_{\v n}^{\v \alpha,\beta}$ taking
$\alpha_j=\beta+c_j\sqrt{\beta}$, the substitution $x\to
\frac{x+\sqrt{\beta}}{2\sqrt{\beta}}$ and letting $\beta\to
+\infty$ after multiplying with some constant depending on $\v n$
and $\beta$.  It seems that we can't find explicit expressions for
the multiple Hermite polynomials in terms of a $\mathcal{M}_{q,\v
n}^{p;r}$ series or a Kamp\'e de F\'eri\'et series using such
arguments.
\end{itemize}

\section{Conclusion}
\begin{figure}[t]
\begin{center}
\includegraphics[scale=0.7]{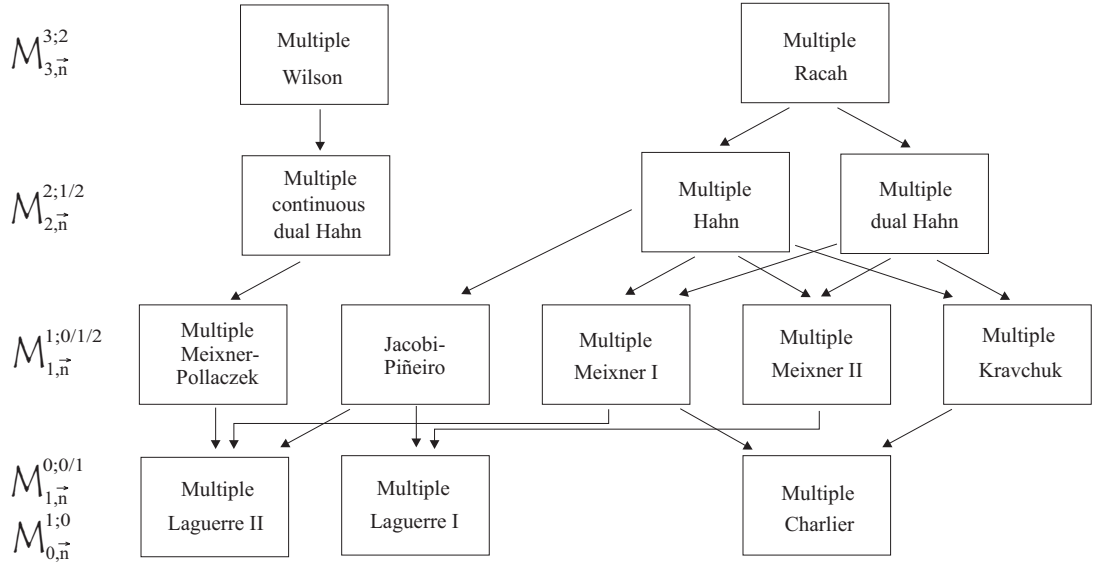}
\end{center}
\caption{\label{Askeyscheme}{\em The (incomplete) multiple
AT-Askey scheme}}
\end{figure}

In Figure \ref{Askeyscheme} all these and some extra limit
relations are combined.  This scheme is a generalization of the
Askey scheme in the scalar case.  In each of the examples of this
scheme the measures have the same support. Although for multiple
Wilson, multiple Racah, multiple continuous dual Hahn and multiple
dual Hahn it is still an open question, we believe that they are
all examples of AT systems which is the reason we call it the
multiple AT-Askey scheme.

This scheme doesn't contain all the possible examples of multiple
orthogonal polynomials which reduce to the orthogonal polynomials
of the Askey scheme.  In \cite{Coussement} the authors also
mentioned some examples of Angelesco systems (with their
hypergeometric expression). It is also possible to change more
than one parameter in the Wilson weight (maybe with some
correlation) in order to find other examples of multiple Wilson
polynomials. Then it is for example possible to obtain multiple
continuous Hahn polynomials (with positive measures on the real
line) using some limit relations.

\end{document}